\documentclass[12pt,reqno]{amsart}
\usepackage{enumerate, latexsym, amsmath, amsfonts, amssymb, amsthm, color}
\def\Z{\Bbb Z}
\def\N{\Bbb N}

\def\l{\left}
\def\r{\right}

\def\t{\text}
\def\f{\frac}

\def\ls{\leqslant}
\def\gs{\geqslant}

\def\Proof{\noindent{\it Proof}}

\theoremstyle{plain}
\newtheorem{theorem}{Theorem}

\newtheorem{lemma}{Lemma}

\theoremstyle{definition}

\theoremstyle{remark}
\newtheorem{remark}{Remark}

 \vspace{4mm}

\begin{document}
 \baselineskip=16pt
\hbox{Preprint, {\tt arXiv:1210.1562}}
\medskip

\title
[On irreducible polynomials over finite fields]
{On irreducible polynomials over finite fields}

\author
[Zhi-Wei Sun] {Zhi-Wei Sun}

\thanks{Supported by the National Natural Science Foundation (grant 11171140)
 of China and the the PAPD of Jiangsu Higher
Education Institutions}

\address {Department of Mathematics, Nanjing
University, Nanjing 210093, People's Republic of China}
\email{zwsun@nju.edu.cn}


 \begin{abstract} For $n=1,2,3,\ldots$ let $N_n(q)$ denote the number of monic irreducible polynomials over the finite field $\Bbb F_q$.
 We show that the sequence $(\sqrt[n]{N_n(q)})_{n>e^{3+7/(q-1)^2}}$ is strictly increasing
 and the sequence $(\sqrt[n+1]{N_{n+1}(q)}/\sqrt[n]{N_n(q)})_{n\gs5.835\times10^{14}}$ is strictly decreasing.
 We also prove that if $q\gs9$ then $(N_{n+1}(q)/N_n(q))_{n\gs1}$ is strictly increasing.
\end{abstract}

\maketitle

\section{Introduction}
\setcounter{lemma}{0}
\setcounter{theorem}{0}
\setcounter{corollary}{0}
\setcounter{remark}{0}
\setcounter{equation}{0}
\setcounter{conjecture}{0}

 Let $p_n$ denote the $n$-th prime for each $n\in\Z^+=\{1,2,3,\ldots\}$. In number theory, Firoozbakht's conjecture asserts that
$$\root n\of{p_n}>\root{n+1}\of {p_{n+1}}\quad\t{for all}\ n=1,2,3,\ldots,$$
i.e., the sequence $(\root{n}\of {p_n})_{n\gs1}$ is strictly
decreasing (cf. \cite[p.\,185]{R}). Though this remains unsolved, recently the author \cite{S} was able to show that
$(\root{n}\of {S_n})_{n\gs2}$ is strictly decreasing
and moreover the sequence $(\root{n+1}\of{S_{n+1}}/\root{n}\of
{S_n})_{n\gs5}$ is strictly increasing, where $S_n$ is the sum of the first $n$ primes.

Let $q>1$ be a prime power and let $\Bbb F_q$ be the finite field of order $q$.
For each $n\in\Z^+$ we use $N_n(q)$ to denote the number of monic irreducible polynomials over $\Bbb F_q$.

In this paper we establish the following new result.

\begin{theorem}\label{Th1.1} Let $q>1$ be any prime power. Then $(\sqrt[n]{N_n(q)})_{n>e^{3+7/(q-1)^2}}$
is strictly increasing, and the sequence
$$(\sqrt[n+1]{N_{n+1}(q)}/\sqrt[n]{N_n(q)})_{n\gs5.835\times10^{14}}$$ is strictly decreasing.
Also,
$(N_{n+1}(q)/N_n(q))_{n\gs1}$ is strictly increasing if $q\gs9$, and
$(N_{n+1}(q)/N_n(q))_{n\gs19}$ is strictly increasing if $q<9$.
\end{theorem}
\begin{remark}\label{Rem1.1} Our computation suggests that $(\sqrt[n]{N_n(q)})_{n\gs2}$ is strictly increasing
and the sequence $(\sqrt[n+1]{N_{n+1}(q)}/\sqrt[n]{N_n(q)})_{n\gs n_0(q)}$ is strictly decreasing, where $n_0(2)=14$, $n_0(3)=8$,
$n_0(4)=n_0(5)=6$ and $n_0(q)=4$ for $q>5$.
\end{remark}

\section{Proof of Theorem 1.1}
\setcounter{lemma}{0}
\setcounter{theorem}{0}
\setcounter{corollary}{0}
\setcounter{remark}{0}
\setcounter{equation}{0}
\setcounter{conjecture}{0}

\begin{lemma}\label{Lem2.1} For any integer $n>1$ we have
\begin{equation}\label{2.1}|nN_n(q)-q^n|<\f{q^{n/p(n)+1}}{q-1},
\end{equation}
where $p(n)$ denotes the least prime divisor of $n$.
\end{lemma}
\Proof. It is well known that
$$N_n(q)=\f1n\sum_{d\mid n}\mu(d)q^{n/d}$$
(see, e.g., \cite[p.\,84]{IR}), where $\mu$ denotes the M\"obius function. So
\begin{align*}|nN_n(q)-q^n|=&\bigg|\sum_{d\mid n\atop d\gs p(n)}\mu(d)q^{n/d}\bigg|
\\\ls&\sum_{d=p(n)}^nq^{n/d}\ls\sum_{k=1}^{n/p(n)}q^k=q\f{q^{n/p(n)}-1}{q-1}<\f{q^{n/p(n)+1}}{q-1}.
\end{align*}
This concludes the proof. \qed

\begin{lemma}\label{Lem2.2} Let $n>1$ be an integer and set $L_n(q)=(q-1)q^{n-n/p(n)-1}$. Then
\begin{equation}\label{2.2}\l|\log N_n(q)-\log\f{q^n}n\r|<\f2{L_n(q)}.
\end{equation}
\end{lemma}
\Proof. Clearly $N_2(2)=1$ and hence (\ref{2.2}) holds in the case $q=n=2$.

Below we assume that $q>2$ or $n>2$. If $n>2$ then $n-n/p(n)\gs2$. So $L_n(q)\gs2$.
Write $nN_n(q)=q^n(1+r_n(q))$. By Lemma 2.1, $|r_n(q)|<1/L_n(q)\ls1/2$.
If $r_n(q)\gs0$, then
$$0\ls\log N_n(q)-\log\f{q^n}n=\log(1+r_n(q))\ls r_n(q)<\f1{L_n(q)}.$$
As $\log(1-x)>-2x$ for $x\in(0,1/2)$, when $r_n(q)<0$ we have
$$0>\log N_n(q)-\log\f{q^n}n=\log(1-|r_n(q)|)>-2|r_n(q)|>\f{-2}{L_n(q)}.$$
Therefore (\ref{2.2}) always holds. \qed

\begin{lemma}\label{Lem2.3} Let $n>4$ be an integer. We have
\begin{equation}\label{2.3}L_n(q)\gs q-1+\f{n-2}2(q-1)^2+\f{(n-2)(n-4)}8(q-1)^3>\f{(n-1)^2}8.
\end{equation}
If $n\gs6$, then
\begin{equation}\begin{aligned}L_n(q)\gs &q-1+\f{n-2}2(q-1)^2+\f{n^2-6n+8}{8}(q-1)^3
\\&+\f{n^3-12n^2+44n-48}{48}(q-1)^4.\end{aligned}\end{equation}
\end{lemma}
\Proof. As $n-n/p(n)-1\gs2$, we have
\begin{align*}&\f{L_n(q)}{q-1}=(1+(q-1))^{n-n/p(n)-1}
\\\gs&1+\l(n-\f n{p(n)}-1\r)(q-1)+\l(n-\f n{p(n)}-1\r)\l(n-\f n{p(n)}-2\r)\f{(q-1)^2}2
\\\gs&1+\l(\f n2-1\r)(q-1)+\l(\f n2-1\r)\l(\f n2-2\r)\f{(q-1)^2}2
\end{align*}
and hence
\begin{align*}L_n(q)\gs& q-1+\f{n-2}2(q-1)^2+\f{(n-2)(n-4)}8(q-1)^3
\\\gs&1+\f{n-2}2+\f{n^2-6n+8}8>\f{(n-1)^2}8.
\end{align*}
Similarly, if $n\gs6$ then
\begin{align*}&\f{L_n(q)}{q-1}=(1+(q-1))^{n-n/p(n)-1}
\\\gs&1+\l(n-\f n{p(n)}-1\r)(q-1)+\l(n-\f n{p(n)}-1\r)\l(n-\f n{p(n)}-2\r)\f{(q-1)^2}2
\\&+\l(n-\f n{p(n)}-1\r)\l(n-\f n{p(n)}-2\r)\l(n-\f n{p(n)}-3\r)\f{(q-1)^3}{3!}
\\\gs&1+\l(\f n2-1\r)(q-1)+\l(\f n2-1\r)\l(\f n2-2\r)\f{(q-1)^2}2
\\&+\l(\f n2-1\r)\l(\f n2-2\r)\l(\f n2-3\r)\f{(q-1)^3}6
\end{align*}
which gives the desired (2.4). \qed

\medskip
\noindent{\it Proof of Theorem 1.1}. By Lemma 2.2, for $n=2,3,\ldots$ we may write $\log N_n(q)-\log(q^n/n)=c_n/L_n(q)$ with $|c_n|<2$.

(i) Now let $n>\lfloor e^3\rfloor=20$. Observe that
\begin{align*}&\f{\log N_{n+1}(q)}{n+1}-\f{\log N_n(q)}n
\\=&\f1{n+1}\l(\log\f{q^{n+1}}{n+1}+\f{c_{n+1}}{L_{n+1}(q)}\r)-\f1n\l(\log\f{q^n}n+\f{c_n}{L_n(q)}\r)
\\=&\f{\log n}n-\f{\log(n+1)}{n+1}+\f{c_{n+1}}{(n+1)L_{n+1}(q)}-\f{c_n}{nL_n(q)}
\\>&\f{\log n}{n(n+1)}-\f{\log(1+1/n)}{n+1}-\f2{(n+1)L_{n+1}(q)}-\f2{nL_n(q)}
\\>&\f{\log n}{n(n+1)}-\f1{n(n+1)}-\f2{(n+1)L_{n+1}(q)}-\f2{nL_n(q)}.
\end{align*}
In view of Lemma 2.3,
$$nL_n(q)>\f{n(n-2)}2(q-1)^2>\f 37n(n+1)(q-1)^2$$
and
$$(n+1)L_{n+1}(q)>\f{(n+1)(n-1)}2(q-1)^2\gs\f{10}{21}n(n+1)(q-1)^2.$$
Therefore
\begin{align*}\log\f{\sqrt[n+1]{N_{n+1}(q)}}{\sqrt[n]{N_n(q)}}>&\f{\log n-1-2(7/3+21/10)/(q-1)^2}{n(n+1)}
\\>&\f{\log n-1-9/(q-1)^2}{n(n+1)}.
\end{align*}
If $n>e^{3+7/(q-1)^2}$, then $\log n>1+9/(q-1)^2$ and hence
$\sqrt[n+1]{N_{n+1}(q)}>\sqrt[n]{N_n(q)}$.

(ii) Now we fix an integer $n\gs 5.835\times 10^{14}$ and set
$$\Delta_n(q):=\f2{n+1}\log N_{n+1}(q)-\f{\log N_{n}(q)}{n}-\f{\log N_{n+2}(q)}{n+2}.$$
Observe that
\begin{align*}\Delta_n(q)=&\f2{n+1}\l(\log\f{q^{n+1}}{n+1}+\f{c_{n+1}}{L_{n+1}(q)}\r)
-\f1{n}\l(\log\f{q^{n}}{n}+\f{c_{n}}{L_{n}(q)}\r)
\\&-\f1{n+2}\l(\log\f{q^{n+2}}{n+2}+\f{c_{n+2}}{L_{n+2}(q)}\r)
\\=&\f{\log n}{n}+\f{\log(n+2)}{n+2}-\f2{n+1}\log (n+1)
\\&+\f{2c_{n+1}}{(n+1)L_{n+1}(q)}-\f{c_{n}}{nL_{n}(q)}-\f{c_{n+2}}{(n+2)L_{n+2}(q)}
\end{align*}
and hence
\begin{align*}\Delta_n(q)>&\l(\f1{n}+\f1{n+2}-\f2{n+1}\r)\log n+\f{\log(1+2/n)}{n+2}-\f2{n+1}\log\l(1+\f1n\r)
\\&-\f{4}{(n+1)L_{n+1}(q)}-\f{2}{nL_{n}(q)}-\f2{(n+2)L_{n+2}(q)}
\\>&\f{2\log n}{n(n+1)(n+2)}+\f1{n+2}\l(\f2n-\f2{n^2}\r)-\f2{n+1}\cdot\f1n
\\&-\f{4}{(n+1)L_{n+1}(q)}-\f{2}{nL_{n}(q)}-\f2{(n+2)L_{n+2}(q)}.
\end{align*}
Therefore
\begin{equation}\label{2.5}\begin{aligned}\f{\Delta_n(q)}2>&\f{\log n-1}{n(n+1)(n+2)}-\f1{n^2(n+2)}
\\&-\f{2}{(n+1)L_{n+1}(q)}-\f{1}{nL_{n}(q)}-\f1{(n+2)L_{n+2}(q)}.
\end{aligned}\end{equation}

As $n>58\times 10^{13}$, by Lemma 2.3 we have
\begin{gather*}\f1{n^2(n+2)}<\f{1+10^{-14}}{n(n+1)(n+2)},
\\\f1{nL_n(q)}<\f 8{n(n-1)^2}<\f{8+10^{-13}}{n(n+1)(n+2)},
\\\f1{(n+1)L_{n+1}(q)}<\f 8{(n+1)n^2}<\f{8+10^{-13}}{n(n+1)(n+2)},
\\\f1{(n+2)L_{n+2}(q)}<\f{8}{(n+2)(n+1)^2}<\f{8}{n(n+1)(n+2)}.
\end{gather*}
Combining these with (2.5) we obtain
\begin{align*}\f{n(n+1)(n+2)}2\Delta_n(q)>&\log n-1-(1+10^{-14})-3(8+10^{-13})-8
\\=&\log n-(34+31\times10^{-14})>0
\end{align*}
since $n\gs 5.835\times 10^{14}>e^{34+31/10^{14}}$.
Thus
$$\log\f{\sqrt[n+1]{N_{n+1}(q)}}{\sqrt[n]{N_n(q)}}>\log\f{\sqrt[n+2]{N_{n+2}(q)}}{\sqrt[n+1]{N_{n+1}(q)}}$$
as desired.

(iii) Observe that
\begin{align*}&2\log N_{n+1}(q)-\log N_n(q)-\log N_{n+2}(q)
\\=&2\l(\log\f{q^{n+1}}{n+1}+\f{c_{n+1}}{L_{n+1}(q)}\r)-\l(\log\f{q^n}n+\f{c_n}{L_n(q)}\r)-\l(\log\f{q^{n+2}}{n+2}+\f{c_{n+2}}{L_{n+2}(q)}\r)
\\=&\log\l(1-\f1{(n+1)^2}\r)+\f{2c_{n+1}}{L_{n+1}(q)}-\f{c_n}{L_n(q)}-\f{c_{n+2}}{L_{n+2}(q)}
\\<&-\f1{(n+1)^2}+\f{4}{L_{n+1}(q)}+\f{2}{L_n(q)}+\f{2}{L_{n+2}(q)}
\end{align*}
If $n\gs6$ and $q\gs9$, then by (2.4) we have
$$L_n(q)\gs8+\f{n-2}2 8^2+\f{n^2-6n+8}88^3+\f{n^3-12n^2+44n-48}{48}8^4=\f83P(n),$$
where $$P(n)=32n^3-360n^2+1276n-1365.$$
Thus, when $n\gs6$ and $q\gs9$, we have
\begin{align*}&2\log N_{n+1}(q)-\log N_n(q)-\log N_{n+2}(q)
\\<&-\f1{(n+1)^2}+\f{4\times 3}{8P(n+1)}+\f{2\times 3}{8P(n)}+\f{2\times 3}{8P(n+2)}<0.
\end{align*}
Since
\begin{gather*}N_1(q)=q,\ N_2(q)=\f{q(q-1)}2,
\ N_3(q)=\f{q(q^2-1)}3,\ N_4(q)=\f{q^2(q^2-1)}4,
\\N_5(q)=\f{q(q^4-1)}5,\ N_6(q)=\f{q^6-q^3-q^2+q}6,\ N_7(q)=\f{q^7-q}7,
\end{gather*}
we can easily verify that $N_{n+1}(q)^2<N_n(q)N_{n+2}(q)$ for $n=1,2,3,4,5$.

If $n\gs6$ and $q\ls8$, then by (2.4) we have
\begin{align*}L_n(q)\gs&1+\f{n-2}2 +\f{n^2-6n+8}8+\f{n^3-12n^2+44n-48}{48}
\\=&\f n{48}(n^2-6n+32);
\end{align*}
similarly,
$$L_{n+1}(q)\gs\f{n+1}{48}(n^2-4n+27)\ \ \t{and}\ \ L_{n+2}(q)\gs\f{n+2}{48}(n^2-2n+24).$$
Thus, when $n\gs389$ and $q\ls8$, we have
\begin{align*}&2\log N_{n+1}(q)-\log N_n(q)-\log N_{n+2}(q)
\\<&-\f1{(n+1)^2}+\f{4\times 48}{(n+1)(n^2-4n+27)}
\\&+\f{2\times 48}{n(n^2-6n+32)}+\f{2\times 48}{(n+2)((n^2-2n+24)}
\\<&0.
\end{align*}
If $q\ls8$ and $19\ls n\ls 388$, then it is easy to verify  $N_{n+1}(q)^2<N_n(q)N_{n+2}(q)$ via computer.

In view of the above, we have completed the proof of Theorem 1.1.  \qed

\end{document}